\documentclass[12pt, a4paper]{article}
\usepackage[latin1]{inputenc}
\usepackage[english]{babel}
\usepackage{indentfirst}
\usepackage{amstext}
\usepackage{setspace}
\usepackage{amsfonts}
\usepackage{textcomp}
\usepackage{amssymb}
\usepackage{amscd}
\usepackage{epsf}
\usepackage{graphicx}
\usepackage{epsfig}
\usepackage{MnSymbol}

\newcommand {\dem}{\hskip -0.6cm{\bf Proof:  }}
\newcommand {\fim}{\hfill{$\square$}\vskip 1pc}

\newcommand {\N}{\mathbb{N}}

\newcommand {\F}{\mathbb{F}}
\newcommand {\G}{\mathrm{G}}

\newcommand{\supp}{\text{supp}}

\newtheorem{theorem}{Theorem}[section]
\newtheorem{lema}[theorem]{Lemma}

\newtheorem{definition}[theorem]{Definition}
\newtheorem{proposition}[theorem]{Proposition}
\newtheorem{example}[theorem]{Example}
\newtheorem{obs}[theorem]{Remark}

\DeclareMathOperator{\idmap}{id}
\DeclareMathOperator{\Kspan}{span}
\newcommand{\augmap}{\mathcal{T}}

\begin{document}

\title{Simplicity of partial skew group rings and maximal commutativity}
\maketitle
\begin{center}
%{\Large{\bf Perron-Frobenius operators and representations of the Cuntz-Krieger algebras for infinite matrices}} \\
%\vspace{5mm}
{\large Daniel Gonçalves, Johan \"{O}inert and Danilo Royer}\\
\end{center}  
\vspace{8mm}

\abstract 
Let $R_0$ be a commutative associative ring (not necessarily unital), $G$ a group and $\alpha$ a partial action by ideals that contain local units. We show that $R_0$ is maximal commutative in the partial skew group ring $R_0\rtimes_\alpha G$ if and only if $R_0$ has the ideal intersection property in $R_0\rtimes_\alpha G$. From this we derive a criterion for simplicity of $R_0\rtimes_\alpha G$ in terms of maximal commutativity and $G-$simplicity of $R_0$ and apply this to two examples, namely to partial actions by clopen subsets of a compact set and to give a new proof of the simplicity criterion for Leavitt path algebras. A new proof of the Cuntz-Krieger uniqueness theorem for Leavitt path algebras is also provided.

\doublespace


\onehalfspace

\section{Introduction}

Partial skew group rings arose as a generalization of skew group rings and as an algebraic analogue of C*-partial crossed products (see \cite{Ex}). Much in the same way as skew group rings, partial skew group rings provide a way to construct non-commutative rings, and recently Leavitt path algebras have been realized as partial skew group rings (see \cite{rg1}), indicating that the theory of non-commutative rings may benefit from the theory of partial skew group rings. Still, when compared to the well-established theory of skew group rings, the theory of partial skew group rings is still in its infancy. Actually, to our knowledge, \cite{Ferrero} and \cite{BG} are the only existing papers regarding the ideal structure of partial skew group rings, and \cite{DG} is a recent paper describing simplicity conditions for partial skew group rings of abelian groups. 

Our main goal in this paper is to derive necessary and sufficient conditions for simplicity of partial skew group rings. In general, this is still an open problem, even for skew group rings. In \cite{Oinert1} and \cite{Oinert},  \"{O}inert has attacked this problem for skew group rings $R_0\rtimes_\alpha G$, where either the group $G$, or the ring $R_0$, is abelian. Recently, in \cite{DG}, a criterion for simplicity of partial skew group rings of abelian groups has been described. In our case, we will extend results of \cite{Oinert1} to partial skew group rings $R_0\rtimes_\alpha G$, where $R_0$ is assumed to be commutative and associative (not necessarily unital) and $\alpha$ is a partial action by ideals that contain local units. More specifically, we will show that $R_0\rtimes_\alpha G$ is simple if and only if $R_0$ is $G-$simple and maximal commutative in $R_0\rtimes_\alpha G$. In particular, our results can be applied to Leavitt path algebras, by realizing them as partial skew group rings (see \cite{rg1}), and to partial skew group rings associated with partial topological dynamics.

Our work is organized in the following way: In section 2 we present our main results, preceded by a quick overview of the key concepts involved below. In section 3 we apply the results of section 2 to derive a new proof of the simplicity criterion for Leavitt path algebras, as well as a new proof of the Cuntz-Krieger uniqueness theorem for Leavitt path algebras, and in section 4 we show an application of the results of section 2 to partial topological dynamics, namely to partial actions by clopen subsets of a compact set.

Recall that a partial action of a group $G$ on a set $\Omega$ is a pair $\alpha= (\{D_{t}\}_{t\in \G}, \ \{\alpha_{t}\}_{t\in \G})$, where for each $t\in \G$, $D_{t}$ is a subset of $\Omega$ and $\alpha_{t}:D_{t^{-1}} \rightarrow D_{t}$ is a bijection such that $D_{e} = \Omega$, $\alpha_{e}$ is the identity in $\Omega$, $\alpha_{t}(D_{t^{-1}} \cap D_{s})=D_{t} \cap D_{ts}$ and $\alpha_{t}(\alpha_{s}(x))=\alpha_{ts}(x),$ for all $x \in D_{s^{-1}} \cap D_{s^{-1} t^{-1}}.$ In case $\Omega$ is a ring (algebra) then, for each $t\in G$, the subset $D_t$ should be an ideal and the map $\alpha_t$ should be a ring (algebra) isomorphism. In the topological setting each $D_t$ should be an open set and each $\alpha_t$ a homeomorphism and in the C*-algebra setting each $D_t$ should be a closed ideal and each $\alpha_t$ should be a *-isomorphism.

Associated to a partial action of a group $G$ on a ring $A$, we have the partial skew group ring, $A \rtimes_\alpha G$, which is the set of all finite formal sums $\sum\limits_{t \in G} a_t\delta_t$, where, for each $t \in G$, $a_t \in D_t$ and $\delta_t$ is a symbol. Addition is defined in the usual way and multiplication is determined by $(a_t\delta_t)(b_s\delta_s) = \alpha_t(\alpha_{t^{-1}}(a_t)b_s)\delta_{ts}$. An ideal $I$ of $A$ is said to be \emph{$G-$invariant} if $\alpha_g(I\cap D_{g^{-1}})\subseteq I $ holds for all $g\in G$. If $A$ and $\{0\}$ are the only $G-$invariant ideals of $A$, then $A$ is said to be \emph{$G-$simple}.

For $a=\sum\limits_{t \in G} a_t \delta_t \in A \rtimes_{\alpha} G $, the support of $a$, which we denote by $\supp(a)$, is the finite set 
$\{t \in G \, \, : \, \, a_t \neq 0 \}$, and the cardinality of $\supp(a)$ is denoted by $\# \supp(a)$. For $g\in G$, the projection map into the $g$ coordinate, $P_g:  A \rtimes _{\alpha} G  \rightarrow  A$, is given by $P_g \left(\sum\limits_{t \in G} a_t \delta_t \right)=  a_g$ and the augmentation map $\augmap :R_0\rtimes_\alpha G\rightarrow R_0$  is defined by $\augmap \left(\sum\limits_{t\in G} a_t\delta_t \right)=\sum\limits_{t\in G} a_t$.

Recall also that the centralizer of a nonempty subset $S$ of a ring $R$, which we denote by $C_R(S)$, is the set of all elements of $R$ that commute with each element of $S$. If $C_R(S)=S$ holds, then $S$ is said to be a \emph{maximal commutative subring} of $R$. Note that a maximal commutative subring is necessarily commutative.
Following \cite{OinertLundstrom}, a subring $S$ of a ring $R$ is said to have the \emph{ideal intersection property} in $R$, if $S\cap I \neq \{0\}$
holds for each non-zero ideal $I$ of $R$.

By abuse of notation, the identity element of an arbitrary group $G$ will be denoted by $0$.

\section{Maximal commutativity, the ideal intersection property and simplicity}\label{secao2}

This is the key section of our paper. Throughout it we will assume that $R_0$ is a commutative and associative ring and $\alpha$ is a partial action of a group $G$ on the ring $R_0$ such that all ideals contain local units. Thus, by \cite{Ex}, the partial skew group ring $R_0\rtimes_\alpha G$ is also associative. We begin by showing the relationship between maximal commutativity of $R_0$ and the ideal intersection property of $R_0$ in $R_0\rtimes_\alpha G$.

\begin{theorem}\label{intideals}
Let $R_0$ be a commutative associative ring, $G$ a group and $\alpha=\left(\{R_t\}_{t\in G},\{\alpha_t\}_{t\in G}\right)$ a partial action such that, for each $t\in G$, $R_t$ contains a set of local units. Then $R_0\delta_0$ is maximal commutative in $R_0\rtimes_\alpha G$ if and only if $I\cap R_0\delta_0\neq \{ 0 \}$ for each non-zero ideal $I$ of $R_0\rtimes_\alpha G$.
\end{theorem}

\dem 
First suppose that $R_0\delta_0$ is maximal commutative in $R_0\rtimes_\alpha G$ and let $I$ be a non-zero ideal of $R_0\rtimes_\alpha G$. We will show that $I\cap R_0\delta_0\neq \{ 0 \}$.

Let $x=\sum\limits_{t\in F} x_t\delta_t$ be a non-zero element in $I$ such that $\# \supp(x)$ is minimal among non-zero elements of $I$ and assume that $x_t\neq 0$ for each $t\in F\subseteq G$. Pick an $s\in F$, let $e\in R_{s^{-1}}$ be an unit for $\alpha_{s^{-1}}(x_s)$ and define $y:=x\cdot e\delta_{s^{-1}}\in I$. Next we show that $y\in R_0 \delta_0$, but first notice that $y\neq 0$ and $\#\supp(y)\leq \#\supp(x)$, since $x_s\neq 0$ and
$$y=x\cdot e\delta_{s^{-1}}=x_s\delta_s\cdot e\delta_{s^{-1}}+\sum\limits_{t\in F\setminus \{s\}}x_t \delta_t \cdot e\delta_{s^{-1}}=x_s\delta_0+\sum\limits_{t\in F\setminus \{s\}}x_t \delta_t \cdot e\delta_{s^{-1}}.$$ 

Now, let $a\in R_0$ and $z:=a\delta_0\cdot y-y\cdot a\delta_0\in I$. Notice that  $\#\supp(z)<\#\supp(x)$, since $a\delta_0\cdot x_s\delta_0-x_s\delta_0\cdot a\delta_0=0$, and hence, from the minimality of $\#\supp(x)$, we have that $z=0$. But this implies that $a\delta_0\cdot y=y\cdot a\delta_0$ for all $a\in R_0$ and so, by the maximal commutativity of $R_0\delta_0$, we obtain that $y\in R_0\delta_0$ and $I\cap R_0\delta_0\neq \{0\}$ as desired.

Next we show that if $R_0\delta_0$ is not maximal commutative in $R_0\rtimes_\alpha G$ then there exists a non-zero ideal $J$ of $R_0\rtimes_\alpha G$ such that $J\cap R_0\delta_0=\{0\}$.

So, suppose that $R_0\delta_0$ is not maximal commutative. This means that there exists an element $a=\sum\limits_{t\in F} a_t\delta_t\in R_0\rtimes_\alpha G\setminus R_0\delta_0$ such that 
$a\cdot b\delta_0=b\delta_0\cdot a$ for all $b\in R_0$, which is equivalent to $a_t\delta_t\cdot b\delta_0=b\delta_0\cdot a_t\delta_t$ for all $t\in F$ and $b\in R_0$. Evaluating the multiplications in this last equation we obtain that $\alpha_t(\alpha_{t^{-1}}(a_t)b)\delta_t= ba_t\delta_t,$ for all $t\in F$ and $b\in R_0$ and hence $\alpha_t(\alpha_{t^{-1}}(a_t)b)=ba_t=a_tb$ for all $t\in F$ and $b\in R_0$. 

Now, fix a non-identity $g\in F$ such that $a_g\neq 0$ and let $J$ be the ideal of $R_0\rtimes_\alpha G$ generated by the element $a_g\delta_0-a_g\delta_g$. 

Notice that each element of $J$ is a finite sum of elements of the form $b_t\delta_t(a_g\delta_0-a_g\delta_g)c_r\delta_r$, where $b_t\delta_t, c_r\delta_r\in R_0\rtimes_\alpha G$. Moreover, $J\neq 0$, since if $e$ is a local unit for $a_g$ then $e\delta_0(a_g\delta_0-a_g\delta_g)e\delta_0$ is a non-zero element of $J$.

We will show that $J$ has null intersection with $R_0\delta_0$ by showing that $\augmap(J) = 0$. In order to do so, notice that, for $b_t\delta_t$ and $c_r\delta_r\in R_0\rtimes_\alpha G$, we have that 
\begin{align*}
b_t\delta_t(a_g\delta_0-a_g\delta_g)c_r\delta_r & =b_t\delta_t\cdot a_g\delta_0\cdot c_r\delta_r-b_t\delta_t\cdot a_g\delta_g\cdot c_r\delta_r \\
  & = b_t\delta_t\cdot a_gc_r\delta_r-b_t\delta_t\cdot \alpha_g(\alpha_{g^{-1}}(a_g)c_r)\delta_{gr} \\
	& =b_t\delta_t\cdot a_gc_r\delta_r-b_t\delta_t\cdot a_gc_r\delta_{gr}=d\delta_{tr}-d\delta_{tgr},
\end{align*}
where $d=\alpha_t(\alpha_{t^{-1}}(b_t)a_gc_r),$ and hence $\augmap(J)=0$. Since the restriction of $\augmap$ to $R_0\delta_0$ is injective we conclude that $J\cap R_0\delta_{0}=\{0\}$ as desired.
\fim

The above result generalizes \cite[Theorem 3.5]{Oinert1}.

\begin{obs} Notice that, in the above theorem, the associativity of $R_0\rtimes_\alpha G$ was only used to prove that the ideal intersection property of $R_0$ in $R_0\rtimes_\alpha G$ implies maximal commutativity of $R_0$. 
\end{obs}

We can now prove the simplicity criterion for $R_0\rtimes_\alpha G$,
and thereby generalize \cite[Theorem 6.13]{Oinert1}.

\begin{theorem}\label{t1} Let $R_0$ be a commutative associative ring, $G$ a group and $\alpha=\left(\{R_t\}_{t\in G},\{\alpha_t\}_{t\in G}\right)$ a partial action of $G$ on $R_0$ such that, for each $t\in G$, $R_t$ has a set of local units. Then the partial skew group ring $R_0\rtimes_\alpha G$ is simple if and only if $R_0$ is $G-$simple and $R_0\delta_0$ is maximal commutative in $R_0\rtimes_\alpha G$. 
\end{theorem}

\dem 
Suppose first that $R=R_0\rtimes_\alpha G$ is simple. By Theorem \ref{intideals} $R_0\delta_0$ is maximal commutative. We show below that $R_0$ is $G-$simple. 

Let $I$ be a $G-$invariant non-zero ideal of $R_0$. Define $J$ as the set of finite sums $\sum a_t\delta_t$ such that $a_t\in I\cap R_t$ for all $t\in G$, that is, $J=\left\{\sum a_t\delta_t\in R:a_t\in I\cap R_t, \,\,\, t\in G\right\}.$

Notice that $J$ is a non-zero ideal of $R$. Indeed, if $a_r\delta_r\in R$ and $a_t\in I\cap R_t$ then 
$a_r\delta_r \cdot a_t\delta_t=\alpha_r(\alpha_{r^{-1}}(a_r)a_t)\delta_{rt}$. Since $I$ is $G-$invariant, $\alpha_r(\alpha_{r^{-1}}(a_r)a_t)\in I$ and by the definition of a partial action $\alpha_r(\alpha_{r^{-1}}(a_r)a_t)\in R_{rt}$ so that $a_r\delta_r \cdot a_t\delta_t\in J$. Similarly, $J$ is a right ideal of $R$ and so, by the simplicity of $R$ we obtain that $J=R$. Now notice that, from the definition of $J$, $P_0(J)=I$ and from what was done above, $P_0(J)=P_0(R)=R_0$. So $I=R_0$ and $R_0$ is $G-$simple.

Suppose now that $R_0$ is $G-$simple and that $R_0\delta_0$ is maximal commutative in $R$. Let $I$ be a non-zero ideal of $R$. By Theorem \ref{intideals}, $I\cap R_0\delta_0\neq \{0\}$. Let $J=I\cap R_0\delta_0$ and notice that $P_0(J)$ is a non-zero ideal of $R_0$. Next we show that $P_0(J)$ is $G-$invariant. 

Let $a_t\in P_0(J)\cap R_t$ and pick a unit $e$ for $a_t$ in $R_t$. Since $a_t\delta_0\in J$ we have that $\alpha_{t^{-1}}(e)\delta_{t^{-1}}\cdot a_t\delta_0\cdot e\delta_t = \alpha_{t^{-1}}(a_t)\delta_0$ is in $J$ and hence $\alpha_{t^{-1}}(a_t)\in P_0(J)$ and $P_0(J)$ is $G-$invariant. 

Now, since $R_0$ is $G-$simple we have that $P_0(J)=R_0$ and so $J=R_0\delta_0$. In particular, $R_0\delta_0 \subseteq I$. Take $s\in G$, $a_s\in R_s$ and an arbitrary $a_s\delta_s\in R_0\rtimes_\alpha G$. Then, letting $e$ be a local unit for $a_s$ in $R_s$, we have that $a_s\delta_s=e\delta_0\cdot a_s\delta_s\in I$. This shows that $R_0\rtimes_\alpha G=I$, as desired.
\fim

Inspired by \cite[Example 3.4]{FL}, we provide the following example.

\begin{example}
Let $R_0 = K e_1 \oplus K e_2 \oplus K e_3$, where $K$ is a field and $e_1, e_2, e_3$ are orthogonal central idempotents of $R_0$.
Let $C_4$ be the cyclic group of order $4$ with generator $g$ and define a partial action of $C_4$ on $R_0$ by $\alpha_0 = \idmap_{R_0}$,
\begin{align*}
\alpha_g : K e_2 \oplus K e_3 \to K e_1 \oplus K e_2, \quad \alpha_g(e_2)=e_1 \quad &\text{ and } \quad \alpha_g(e_3)=e_2; \\
\alpha_{g^2} : K e_1 \oplus K e_3 \to K e_1 \oplus K e_3, \quad \alpha_{g^2}(e_1)=e_3 \quad &\text{ and } \quad \alpha_{g^2}(e_3)=e_1;\\
\alpha_{g^3} : K e_1 \oplus K e_2 \to K e_2 \oplus K e_3, \quad \alpha_{g^3}(e_1)=e_2 \quad &\text{ and } \quad \alpha_{g^3}(e_2)=e_3.
\end{align*}
There are exactly six proper (non-zero) ideals of $R_0$, namely
\begin{displaymath}
	K e_1, \quad K e_2, \quad K e_3, \quad K e_1 \oplus K e_2, \quad K e_1 \oplus K e_3 \quad \text{and} \quad K e_2 \oplus K e_3,
\end{displaymath}
none of which is $C_4-$invariant. One easily checks this using the definition of $\alpha$. 
Thus, $R_0$ is $C_4-$simple.
Moreover, a short calculation reveals that $R_0\delta_0$ is maximal commutative in the partial skew group ring $R_0 \rtimes_\alpha C_4$.
By Theorem \ref{t1}, we conclude that $R_0 \rtimes_\alpha C_4$ is simple.
\end{example}

\section{A new proof of the simplicity criterion for Leavitt path algebras}

Recently Leavitt path algebras have been described as partial skew group rings, see \cite{rg1}. More precisely, the Leavitt path algebra associated to a graph $E$ has been realized as a partial skew group ring of a commutative algebra by the free group on the edges and so we can apply the characterization of simplicity given in section \ref{secao2} to Leavitt path algebras. This will lead to a new proof of the simplicity criterion for Leavitt path algebras that rely solely on partial skew group ring theory. The details follow below, after we have recalled some of the key definitions given in \cite{rg1}.

Given a field $K$ and a graph $E=(E^1,E^0, r, s)$, $L_K(E)$ will denote, as usual, the Leavitt path algebra associated to $E$ (see \cite{arandapino, gonroy} for example), $W$ is the set of all finite paths and $W^{\infty}$ the set of all infinite paths in $E$. %Define, for all $\xi=\xi_1\cdots\xi_n\in W$ the range and source of $\xi$ by $r(\xi)=r(\xi_n)$ and $s(\xi)=s(\xi_1)$. Moreover, define the source of an element $\nu\in W^{\infty}$ by $s(\nu)=s(\nu_1)$, and the source an range of each vertex $v$ by $r(v)=v=s(v)$. Define the set 
The partial action takes place on the set $$X=\{\xi\in W:r(\xi) \text{ is a sink }\}\cup \{v\in E^0: v\text{ is a sink }\}\cup W^{\infty}$$ and the group acting is the free group generated by $E^1$, which is denoted by $\F$. 

The exact definition of the partial action is a bit cumbersome but we reproduce it here for completeness. For each $c\in \F$, let $X_c$ be defined as follows:
\begin{itemize}
\item $X_0:=X$, where $0$ is the neutral element of $\F$. 

\item $X_{b^{-1}}:=\{\xi\in X: s(\xi)=r(b)\},$ for all $b\in W$. 

\item $X_a:=\{\xi\in X: \xi_1\xi_2...\xi_{|a|}=a\},$ for all $a\in W$.

\item $X_{ab^{-1}}:=\{\xi\in X: \xi_1\xi_2...\xi_{|a|}=a\}=X_a,$ for $ab^{-1}\in \F$ with $a,b\in W$, $r(a)=r(b)$ and $ab^{-1}$ in its reduced form.

\item $X_c:=\emptyset$, for all other $c \in \F$.

\end{itemize}

Let  $\theta_{0}:X_{0}\rightarrow X_{0}$ be the identity map. For $b\in W$,  $\theta_b:X_{b^{-1}}\rightarrow X_b$ is defined by $\theta_b(\xi)=b\xi$ and $\theta_{b^{-1}}:X_b\rightarrow X_{b^{-1}}$ by
$\theta_{b^{-1}}(\eta)= \eta_{|b|+1}\eta_{|b|+2}...$ if $r(b)$ is not a sink  and $\theta_{b^{-1}}(b)=r(b)$, if $r(b)$ is a sink. Finally, for $a,b\in W$ with $r(a)=r(b)$ and $ab^{-1}$ in reduced form, $\theta_{ab^{-1}}:X_{ba^{-1}}\rightarrow X_{ab^{-1}}$ is defined by $\theta_{ab^{-1}}(\xi)=a\xi_{(|b|+1)}\xi_{(|b|+2)}...$, with inverse 
$\theta_{ba^{-1}}:X_{ab^{-1}}\rightarrow X_{ba^{-1}}$ defined by $\theta_{ba^{-1}} (\eta)=b\eta_{(|a|+1)}\eta_{(|a|+2)}...$ .

Notice that $\{\{X_c\}_{c\in \F}, \{\theta_c\}_{c\in\F}\}$ is a partial action on the set level and so it induces a partial action $\{\{F(X_c)\}_{c\in \F}, \{\alpha_c\}_{c\in\F}\}$, where, for each $c\in \F$, $F(X_c)$ denotes the algebra of all functions from $X_c$ to $K$, and $\alpha_c:F(X_{c^{-1}})\rightarrow F(X_c)$ is defined by $\alpha_c(f)=f\circ \theta_{c^{-1}}$. The skew group ring associated to this partial action is not $L_K(E)$ yet. For this one proceeds in the following way:

For each $c\in \F$, and for each $v\in E^0$, define the characteristic maps $1_c:=\chi_{X_c}$ and $1_v:=\chi_{X_v}$, where $X_v=\{\xi\in X:s(\xi)=v.\}$. Notice that $1_c$ is the unit of $F(X_c)$. Finally, let
$$D_0=\Kspan\{\{1_p:p\in \F\setminus\{0\}\}\cup\{1_v:v\in E^0\}\},$$  (where $\Kspan$ means the $K$-linear span) and, for each $p\in \F\setminus\{0\}$, let $D_p\subseteq F(X_p)$ be defined as $1_p D_0$, that is, $$D_p=\Kspan\{1_p1_q:q\in \F\}.$$
Since $\alpha_p(1_{p^{-1}}1_q)=1_p1_{pq}$ (see \cite{rg1}), consider, for each $p\in \F$, the restriction of $\alpha_p$ to $D_{p^{-1}}$. Notice that $\alpha_{p}:D_{p^{-1}}\rightarrow D_p$ is an isomorphism of $K$-algebras and, furthermore, $\{\{\alpha_p\}_{p\in \F}, \{D_p\}_{p\in \F}\}$ is a partial action. In \cite{rg1} it is shown that the partial skew group ring $D_0\rtimes_\alpha\F$ is isomorphic to the Leavitt path algebra $L_K(E)$.

Recall, see \cite{tomforde}, that a subset $H\subseteq E^0$ is said to be hereditary if for any $e\in E^1$ we have that $s(e)\in H$ implies $r(e)\in H$. A hereditary subset $H \subseteq E_0$ is called saturated if whenever $0< \#s^{-1}(v)<\infty$, then $\{r(e)\in H: e\in E^1 \text{ and } s(e) =v \} \subseteq H$ implies $v\in H$.  In \cite{tomforde} it is proved that $L_K(E)$ is simple if and only if the graph $E$ satisfies condition $(L)$, that is, each closed path in the graph $E$ has an exit, and the only hereditary and saturated subsets of $E^0$ are $E^0$ and $\emptyset$. 
From now until the end of this section we will focus on the proof of the above simplicity criterion for $D_0\rtimes_\alpha\F$ via Theorem \ref{t1}, thus giving a new proof of the simplicity criterion for Leavitt path algebras. 
On the way, we will obtain some useful results that we will also use, together with Theorem \ref{intideals}, to give a new proof of the Cuntz-Krieger uniqueness theorem for Leavitt path algebras.

%\begin{definition} Let $E$ be a directed graph, and $K$ be a field. The Leavitt path algebra of $E$ with coefficients in $K$, denoted by $L_K(E)$, is the universal $K$-algebra generated by a set $\{v: v \in E^0\}$, of pairwise orthogonal idempotents, together with a set $\{e, e^* : e \in E^1\}$ of elements satisfying:
%\begin{enumerate}
%\item $s(e)e = er(e) = e$ for all $e\in E^1$,
%\item $r(e)e^* = e^*s(e) = e^*$ for all $e\in E^1$,
%\item for all $e,f\in E^1$, $e^*f =0$ if $e\neq f$ and $e^*e=r(e)$.
%\item $v =\sum\limits_{e\in E^1:s(e)=v}ee^*$ for every vertex $v$ with $0 < \#\{e: s(e) = v\} <\infty.$
%\end{enumerate}
%\end{definition}

%By \cite[3.2]{rg1} and \cite[3.3]{rg1} we obtain the following Theorem:

%\begin{theorem}\label{t4} The K-algebras $L_K(E)$ and $D_0\rtimes_\alpha\F$ are K-isomorphic, via a K-isomorphism $\varphi$, from $L_K(E)$ on $D_0\rtimes_\alpha\F,$ defined on the generators of $L_K(E)$ by $\varphi(e)=1_e\delta_e$, $\varphi(e^*)=1_{e^{-1}}\delta_{e^{-1}}$, for all $e\in E^1$, and $\varphi(v)=1_v\delta_0$, for all $v\in E^0$.
%\end{theorem}

\begin{proposition}\label{t2}
The set $D_0\delta_0$ is maximal commutative in $D_0\rtimes_\alpha \F$ if and only if the graph $E$ satisfies condition $(L)$.
\end{proposition}

\dem Suppose first that $E$ satisfies condition $(L)$. We will show that $D_0\delta_0$ is maximal commutative by contradiction. For this, suppose that there exists an element $a_t\in D_t$, with $t\neq 0$ and $a_t\neq 0$, such that $a_t\delta_t\cdot a_0\delta_0=a_0\delta_0\cdot a_t\delta_t$ for each $a_0\in D_0$, that is, such that \begin{equation}\label{eq1} \alpha_t(\alpha_{t^{-1}}(a_t)a_0)=a_ta_0 \end{equation} for all $a_0\in D_0$. 

Notice that $a_t \neq 0$ implies that either $t\in W$ or $t=r^{-1}$, with $r\in W$, or  $t=ab^{-1}$, where $a,b\in W$.
Furthermore, if in equation \eqref{eq1} we take $a_0=1_{t^{-1}}$ we obtain that $a_t=a_t1_{t^{-1}}$ and hence the support of $a_t$ is contained in $D_t \cap D_{t^{-1}}$ and so $t$ must be a closed path. 

Now, taking appropriate functions for $a_0$ in equation \eqref{eq1} and using induction we obtain that, for all $n\in \N$, $a_t=a_t1_{(t^n)^{-1}}$ and $a_t1_{t^n}=a_t$. For example, for $a_0=1_{t^{-1}t^{-1}}$ we obtain that $a_t1_{t^{-1}}=a_t1_{t^{-1}t^{-1}}$ and so $a_t=a_t1_{t^{-1}t^{-1}}$. On the other hand, for $a_0 = 1_t 1_{t^{-1}}$ we get that $\alpha_t(\alpha_{t^{-1}}(a_t) 1_t 1_{t^{-1}})= a_t 1_t 1_{t^{-1}}$ and hence $a_t 1_{tt} = a_t 1_{t^{-1}} = a_t$. 

Before we derive our contradiction, notice that if $\xi\in X_t$ is such that $a_t(\xi)\neq 0$ then, since $a_t\in D_t$, there exists an $m\in \N$ such that for each $\mu\in X_t$ with $\mu_1\cdots \mu_m=\xi_1\cdots \xi_m$ it holds that $a_t(\mu)=a_t(\xi)$. We now separate our argument into three cases.

{\it Case 1: Suppose $t\in W$.}

Since $a_t=a_t1_{t^m}$ then $t^m=\xi_1\cdots\xi_m \cdots \xi_{m|t|}$. Let $s$ be an exit for $t$ and $\mu\in X_t$ be such that $\mu_1 \cdots \mu_{m|t|}\cdots \mu_k=t^mt_1...t_ls$. Then $a_t(\mu)=a_t(\xi)\neq 0$, but $a_t(\mu)=a_t(\mu)1_{t^{m+1}}(\mu)=0$, a contradiction. So $t$ is not an element of $W$.

{\it Case 2: Suppose $t=r^{-1}$, with $r\in W$.}

This case follows as the previous one, by using the equality $a_t=a_t1_{(t^{m})^{-1}}$ instead of $a_t=a_t1_{t^m}$.

{\it Case 3: Suppose $t=ab^{-1}$, where $a,b\in W$.}

We obtain a contradiction by proceeding as in case 1 if $|a|\geq |b|$ and as in case 2 if $|a|<|b|$. The details are left to the reader.

We conclude that there is no $a_t \in D_t$, with $t\neq 0$, such that $a_t\delta_t$ commutes with each element of $D_0\delta_0$ and hence $D_0\delta_0$ is maximal commutative.

Suppose now that $E$ does not satisfy condition $(L)$, that is, there exists a closed path $t=t_1...t_m$ which has no exit. We will show that $1_t\delta_t$ commutes with $D_0\delta_0$ and so $R_0 \delta_0$ is not maximal commutative.

Recall that $D_0=\Kspan\{\{1_p:p\in \F\setminus \{0\}\}\cup\{1_v:v\in E^0\}\}$ and so it is enough to show that $1_t\delta_t$ commutes with $1_v\delta_0$ and with $1_p\delta_0$, for each $v\in E^0$ and $p\in \F\setminus\{0\}$.

Let $v\in E^0$. Then $1_t\delta_t\cdot1_v\delta_0=\alpha_t(\alpha_{t^{-1}}(1_t)1_v)\delta_t=\alpha_t(1_{t^{-1}}1_v)\delta_t$  which, by \cite[Lemma 2.3(2)]{rg1}, is non-zero only if $r(t)=v$, in which case is equal to $1_t\delta_t$. On the other hand, $1_v\delta_0\cdot 1_t\delta_t=1_v1_t\delta_t,$ which is non-zero only if $s(t)=v$, in which case is equal to $1_t\delta_t$. Since $t$ is a closed path it follows that $1_t\delta_t$ commutes with $1_v\delta_0$.

Now let $r\in \F\setminus\{0\}$. Notice that, in order to check that $1_t\delta_t$ commutes with $1_r\delta_0$ it is enough to verify that $\alpha_t(1_{t^{-1}}1_r) = 1_t1_r$, which is equivalent to $1_t1_{tr}=1_t1_r$ (since $\alpha_t(1_{t^{-1}}1_r)=1_t1_{tr}$). As before, we now divide our proof into cases:

{\it Case 1: $r\in W$.} If $r=t^nt_1...t_k$ for some $n\geq 0$ and $1\leq k\leq m$ then, since $t$ has no exit, $X_r = X_t = \{tttt\cdots\}$ and hence   $1_t1_{tr}=1_t=1_t1_r$. If $r \in W$ is not of the above form, then $1_t1_{tr}=0=1_t1_r$.

{\it Case 2: $r=s^{-1}$ with $s\in W$.} Suppose first that $r(s)=r(t)$. Then $X_{s^{-1}}= X_t$, since $t$ is a closed path with no exit, and hence $1_t1_{tr}=1_t1_{ts^{-1}}=1_t=1_t1_{s^{-1}}=1_t1_r$. If $r(s)\neq r(t)$, then $1_{ts^{-1}}=0=1_t1_{s^{-1}}$.

{\it Case 3: $r=ab^{-1}$ with $a,b\in W$ and $r(a)=r(b)$.} Since $1_{tr}=1_{tab^{-1}}=1_{ta}$ and $1_{r}=1_{ab^{-1}}=1_a$ this case reduces to case 1.

{\it Case 4: All other $r \in \F$.} In this case $1_r = 0$ and hence both sides of the equation $\alpha_t(1_{t^{-1}}1_r) = 1_t1_r$ are equal to zero. 

We have proved that $1_t \delta_t$ is in the centralizer of $D_0 \delta_0$ and hence $D_0 \delta_0$ is not maximal commutative, as desired.
\fim

Before we proceed to show the connection between $\F-$simplicity of $D_0$ and the nonexistence of proper hereditary and saturated subsets of $E^0$,
we shall prove two useful lemmas.

\begin{lema}\label{NewLemma}
Let $x_0 \delta_0$ be a non-zero element of $D_0 \delta_0$ and denote by $I$ the principal ideal of $D(X) \rtimes_\alpha \F$ generated by $x_0 \delta_0$. Then there exists a vertex $v \in E^0$ such that $1_v \delta_0 \in I$.
\end{lema}

% ALTERNATIVE 1: LONG VERSION
\dem We can write $x_0$ as a linear combination of characteristic functions;
$x_0 = \sum_{i=1}^n \lambda_i 1_{a_i b_i^{-1}} + \sum_{j=1}^m \beta_j 1_{v_j}$, where
$a_i\in W$ and $b_i \in W\cup \{0\}$ (if $a_i=0$, then $1_{a_i b_i^{-1}} = 1_{b_i} = 1_{r(b_i)}$).
Choose some $v\in E^0$ such that $1_v x_0 \neq 0$.
If $v$ is a sink, then $1_v 1_{a_i b_i^{-1}} = 0$ for each $i$, and then
\begin{equation*}
	0 \neq 1_v x_0 \delta_0 = \sum_{j=1}^n \beta_j 1_v 1_{v_j} \delta_0 = \sum_{j : v_j=v} \beta_j 1_v \delta_0
\label{eq:sinkcalculation}
\end{equation*}
which shows that $1_v \delta_0 \in I$.

Now, suppose that $v$ is not a sink. Let $m=\max\{|a_i| \mid 1 \leq i \leq n\}$.
Recall that we can write $X_v = \bigcupdot_{c\in I} X_c$ where the index set $I$ consists of all
all $c\in W$ such that $s(c)=v$ and $|c|=m$ or $s(c)=v$, $|c|<m$ and $r(c)$ is a sink.
If $1_c 1_{a_i b_i^{-1}} \neq 0$, then $a_i$ is the beginning of $c$, and then
$1_c 1_{a_i b_i^{-1}} = 1_c 1_{a_i} = 1_c$. Moreover, if $1_c 1_{v_j}\neq 0$, then $1_c 1_{v_j} = 1_c$.
Using this, we obtain
\begin{align*}
0 \neq 1_c x_0 \delta_0 = \sum_{i=1}^n \lambda_i 1_c 1_{a_i b_i^{-1}} \delta_0 + \sum_{j=1}^m \beta_j 1_c 1_{v_j} \delta_0 = \\
= \sum_{i : 1_c 1_{a_i b_i^{-1}} \neq 0} \lambda_i 1_c 1_{a_i b_i^{-1}} \delta_0 + \sum_{j : 1_c 1_{v_j} \neq 0} \beta_j 1_c 1_{v_j} \delta_0 = \\
= \sum_{i : 1_c 1_{a_i b_i^{-1}} \neq 0} \lambda_i 1_c \delta_0 + \sum_{j : 1_c 1_{v_j} \neq 0} \beta_j 1_c \delta_0
= \left(\sum_{i : 1_c 1_{a_i b_i^{-1}} \neq 0} \lambda_i + \sum_{j : 1_c 1_{v_j} \neq 0} \beta_j \right) 1_c \delta_0.
\end{align*}
which shows that $1_c \delta_0 \in I \setminus \{0\}$.
Note that $1_{r(c)} \delta_0 = 1_{c^{-1}} \delta_0 = 1_{c^{-1}} \delta_{c^{-1}} \cdot 1_c \delta_0 \cdot 1_c \delta_c$.
Using that $I$ is an ideal, we conclude that $1_{r(c)} \delta_0 \in I$ which proves the lemma.
\fim

%% ALTERNATIVE 2: SHORT VERSION
%\dem The proof is similar to the proof of \cite[Proposition 4.3]{rg1}.\fim

\begin{lema}\label{p1} Let $I$ be an $\F-$invariant ideal of $D_0$. Then, the set $Z=\{v\in E^0: 1_v\in I\}$ is hereditary and saturated.
\end{lema}

\dem Let $e\in E^1$ be such that $s(e)\in Z$. Then $1_e=1_{s(e)}1_e\in I\cap D_e$ and, by the $\F-$invariance of $I$, $\alpha_{e^{-1}}(1_e)=1_{e^{-1}}=1_{r(e)}\in I$, so that $r(e)\in Z$.

Now, let $v\in E^0$ be such that $0<\# s^{-1}(v)<\infty$ and $r(e)\in Z$ for each $e\in s^{-1}(v)$. Notice that $1_{r(e)}=1_{e^{-1}}$ and so, since $I$ is $\F-$invariant, we have that $1_e=\alpha_e(1_{e^{-1}})\in I$. This implies that $1_v=\sum\limits_{e\in s^{-1}(v)}1_e\in I$ and hence $v\in Z$ as desired.
\fim

%The theorem above, together with the proposition below, give us the necessity of the simplicity criterion for $D_0\rtimes_\alpha \F$.
The following proposition gives us a characterization of $\F-$simplicity of $D_0$.

\begin{proposition}\label{t3} The algebra $D_0$ is $\F-$simple if and only if the only saturated and hereditary subsets of $E^0$ are $E^0$ and $\emptyset$.
\end{proposition}

\dem Suppose first that $D_0$ is $\F-$simple.
Let $F$ be a nonempty saturated and hereditary subset of $E^0$.
We need to show that $F=E^0$.

Consider the ideal $I$ generated by $\{1_v\delta_0:v\in F\}$ in $D_0\rtimes_\alpha \F$, that is, $I$ is the linear span of all the elements of the form $a_r\delta_r1_v\delta_0b_s\delta_s$, with $v\in F$, $a_r\in D_r$, $b_s\in D_s$ and $r,s\in \F$. Let $J=P_0(D_0\delta_0\cap I)$ and notice that $J$ is a non-zero $\F-$invariant ideal of $D_0$ ($J$ is $\F-$invariant since if $a_t\in J\cap D_t$, then $a_t\delta_0\in I$, so $\alpha_{t^{-1}}(a_t)\delta_0=1_{t^{-1}}\delta_{t^{-1}}\cdot a_t\delta_0\cdot 1_t\delta_t\in I$ and hence $\alpha_{t^{-1}}(a_t)\in J$). Now, since $D_0$ is $\F-$simple we have that $J=D_0$ and, in particular, $1_u\in J$ for each $u\in E^0$. This means that for each $u\in E^0$, $1_u\delta_0\in I$, and so we can write
$$1_u\delta_0=\sum\limits_{t} x_t\delta_t\cdot 1_{v_t}\delta_0\cdot y_{t^{-1}}\delta_{t^{-1}}=\sum\limits_{t} \alpha_t(\alpha_t^{-1}(x_t)1_{v_t}y_{t^{-1}})\delta_0,$$
where the sum above is a finite sum and $v_t\in F$ for each $t$. Multiplying the above equation by $1_u\delta_0$, we obtain
$$1_u\delta_0=\sum\limits_{t\in T}1_u\alpha_t(\alpha_t^{-1}(x_t)1_{v_t}y_{t^{-1}})\delta_0,$$
where $$T:=\{ t\in \F : 1_u\alpha_t(\alpha_t^{-1}(x_t)1_{v_t}y_{t^{-1}})\neq 0\}.$$
In particular, since $1_u\alpha_t(\alpha_t^{-1}(x_t)1_{v_t}y_{t^{-1}})\neq 0$ for each $t\in T$, we have that $1_u1_t\neq 0$ and $1_{v_t}1_{t^{-1}}\neq 0$ for all $t\in T$.

Our aim is to show that each $u\in E^0$ belongs to $F$. So, let $u\in E^0$. If $u=r(b)$ for some path $b$ and $s(b)\in F$ then $u\in F$, since $F$ is hereditary. Moreover, if $0<\# s^{-1}(u)<\infty$ and $r(e)\in F$ for each $e\in s^{-1}(u)$ then $u\in F$, since $F$ is saturated. So, we are left with the cases when there is no path $b$ with $s(b)\in F$ and $r(b)=u$ and either $s^{-1}(u)=\emptyset$, $\# s^{-1}(u)=\infty$, or $0<\# s^{-1}(u)<\infty$ but $r(e)\notin F$ for some $e\in s^{-1}(u)$. We handle these three cases below.

{\it Case 1: $s^{-1}(u)=\emptyset$, and there is no path $b$ with $s(b)\in F$ and $r(b)=u$.}

First notice that since there is no $b\in W$ such that $s(b)\in F$ and $r(b)=u$ then, for each $b\in W$, it holds that either $1_u1_{b^{-1}}=0$ or $1_{v}1_b=0$ for each $v\in F$. Then, by the statement right after the definition of $T$, we obtain that there is no $t\in T$ of the form $t=b^{-1}$ (with $b\in W$). Now, for $t$ of the form $t=ab^{-1}\in \F$, with $a\in W$ and $b\in W\cup\{0\}$, we have that $1_u1_t=0$, since $s(a)\neq u$, and hence $t=ab^{-1}\notin T$. We conclude that $T=\{0\}$, and so
$1_u=1_ux_01_{v_0}y_0$ and it follows that $u=v_0\in F$.

{\it Case 2: $\# s^{-1}(u)=\infty$, and there is no path $b$ with $s(b)\in F$ and $r(b)=u$.}

Here, as in case 1, there is no $t\in T$ of the form $t=b^{-1}$ with $b\in W$. Suppose that $0\notin T$.
Then each $t\in T$ is of the form $t=ab^{-1}$, with $a\in W$ and $b\in W\cup \{0\}$. Since $\# s^{-1}(u)=\infty$, there is an element $\xi\in X$ with $s(\xi)=u$ and $s(\xi)\neq s(a)$ for each $ab^{-1}\in T$. Notice that $1_t(\xi)=0$ for all $t\in T$ and so
 $$1=1_u(\xi)=\sum\limits_{t\in T}1_u\alpha_t(\alpha_t^{-1}(x_t)1_{v_t}y_{t^{-1}})(\xi)=0,$$ which is a contradiction. So $0\in T$ and 
 $1_ux_01_{v_0}y_0\neq 0$, which implies that $u=v_0\in F$.

{\it Case 3: $0<\# s^{-1}(u)<\infty$, and there is no path $b$ with $s(b)\in F$ and $r(b)=u$, and there is an edge $e\in s^{-1}(u)$ such that $r(e)\notin F$.}

Again, as in case 1, there is no $t\in T$ of the form $t=b^{-1}$ with $b\in W$. Suppose, as in case 2, that $0\notin T$. Then, as before, each $t\in T$ is of the form $t=ab^{-1}$, with $a\in W$ and $b\in W\cup \{0\}$.

Now, for each $t\in T$, let $c_t=1_u\alpha_t(\alpha_t^{-1}(x_t)1_{v_t}y_{t^{-1}})$. Since, for each $t=ab^{-1} \in T$, it holds that $1_u1_t\neq 0$ and $1_{v_t}1_{t^{-1}}\neq 0$, we have that $s(a)=u$ and $s(b)=v_t\in F$. The heredity of $F$ now implies that $r(b)\in F$ and since $r(a)=r(b)$ we have that $r(a)\in F$. So, we obtain that $$1_u=\sum\limits_{t\in T}c_t=\sum\limits_{ab^{-1}\in T}c_{ab^{-1}},$$ where $u=s(a)$ and $r(a)\in F$ for all $ab^{-1}\in T$.

Let $z=z_1...z_m$ be a path of maximum length such that $|z|\leq \max\{|a|:ab^{-1}\in T\}$ with $s(z)=u$ and $r(z_i)\notin F$ for each $i\in \{1,...,m\}$. By the hypothesis, such a $z$ exists. Then
multiplying the equation $1_u=\sum\limits_{ab^{-1}\in T}c_{ab^{-1}}$ by $1_z$ we obtain

$$1_z=\sum\limits_{ab^{-1}\in T:|z|<|a|,a_1...a_m=z}c_{ab^{-1}}.$$

Since the sum on the right side is finite, we have that $0<\# s^{-1}(r(z))<\infty$. By the maximality of $|z|$, there is no edge $e\in s^{-1}(r(z))$ such that $r(e)\notin F$. Then, $r(e)\in F$ for all $e\in s^{-1}(r(z))$ and, since $F$ is saturated, we obtain that $r(z)\in F$, a contradiction (since $r(z)=r(z_m)\notin F$).

We conclude that $0\in T$ and, as in case 2, it follows that $u\in F$ as desired.

Suppose now, that the only saturated and hereditary subsets of $E^0$ are $E^0$ and $\emptyset$.
Let $I$ be a non-zero $\F-$invariant ideal of $D_0$. We need to show that $I=D_0$.

Let $J$ be the (non-zero) ideal of $D_0\rtimes_\alpha \F$ consisting of all finite sums $\sum\limits a_t\delta_t$, with $a_t\in D_t\cap I$ ($J$ is an ideal since $I$ is $\F-$invariant) and let $Z=\{v\in E^0:1_v\in I\}$. By Lemma \ref{NewLemma}, there is some $v\in E^0$ such that $1_v\delta_0\in J$, so that $1_v\in I$ (since $J\cap D_0\delta_0=I\delta_0$) and hence $Z$ is nonempty. By Lemma \ref{p1}, $Z$ is hereditary and saturated, and therefore $Z=E^0$. Thus, $1_v\in I$ for each $v\in E^0$ and hence $I=D_0$, as desired.
\fim

%The above result concludes the necessity part of the simplicity criterion for Leavitt path algebras. The sufficient part will follows readily after we prove the following proposition.
Propositions \ref{t2} and \ref{t3} above, enable us to translate the language of Leavitt path algebras into the language of partial skew group rings, and vice versa.
Using this, we shall now give a new proof of the simplicity criterion for Leavitt path algebras.

\begin{theorem} The partial skew group ring $D_0\rtimes_\alpha \F$ is simple if and only if the graph $E$ satisfies condition (L) and the only hereditary and saturated subsets of $E^0$ are $\emptyset$ and $E^0$.
\end{theorem}

\dem
By combining the results from Theorem \ref{t1}, Proposition \ref{t2} and Proposition \ref{t3}, the desired conclusion follows.
\fim

%\dem The necessity part of the theorem was already done, since if $D_0\rtimes_\alpha \F$ is simple then, by Theorem \ref{t1}, $D_0\delta_0$ is maximal commutative and $\F-$simple which, by propositions \ref{t2} and \ref{t3}, means that $E$ satisfies condition $(L)$ and the only hereditary saturated subsets of $E^0$ are $\emptyset$ and $E^0$.
%
%Conversely, suppose that $E$ satisfy condition $(L)$ and the only hereditary saturated subsets of $E^0$ are $\emptyset$ and $E^0$. By proposition \ref{t2} we have that $D_0 \delta_0$ is maximal commutative and so, by theorem \ref{t1}, all we need to do is show that $D_0$ is $\F-$simple.
%
%So, let $I$ be a non-zero $\F-$invariant ideal of $D_0$. Take $J$ as the ideal of $D_0\rtimes_\alpha \F$ consisting of all finite sums $\sum\limits a_t\delta_t$, with $a_t\in D_t\cap I$ ($J$ is an ideal since $I$ is $\F-$invariant) and let $Z=\{v\in E^0:1_v\in I\}$. By \cite[4.3]{rg1}, there is a $v\in E^0$ such that $1_v\delta_0\in J$, so that $1_v\in I$ (since $J\cap D_0\delta_0=I\delta_0$) and hence $Z$ is a nonempty. By Lemma \ref{p1} $Z$ is hereditary saturated and it follows that $Z=E^0$. So $1_v\in I$ for each $v\in E^0$ and hence $I=D_0$ and $D_0$ is $\F-$simple as desired.
%\fim

We end this section by providing an alternative proof of the Cuntz-Krieger uniqueness theorem for Leavitt path algebras (cf. \cite{rg1} and \cite{tomforde}).

\begin{theorem}[Cuntz-Krieger uniqueness theorem]
Let $E$ be a graph that satisfies condition (L).
If $\phi : D_0\rtimes_\alpha \F \to B$ is a $K$-algebra homomorphism
such that $\phi(1_v \delta_0)\neq 0$ for each $v\in E^0$, then $\phi$ is injective.
\end{theorem}

\dem
Suppose that $E$ satisfies condition (L) and that $\phi(1_v \delta_0)\neq 0$ for each $v\in E^0$.
Let $I$ denote the ideal $\ker(\phi)$. Seeking a contradiction, suppose that $I \neq \{0\}$.
Proposition \ref{t2} and Theorem \ref{intideals} now yield $D_0 \delta_0 \cap I \neq \{0\}$.
Let $x_0 \delta_0 \in D_0 \delta_0 \cap I$ be a non-zero element.
By Lemma \ref{NewLemma}, there is some $v\in E^0$ such that $1_v \delta_0 \in I = \ker(\phi)$,
but this is a contradiction. Hence $\ker(\phi)=\{0\}$.
\fim

%************************************Section Topological Dynamics **************************
\section{Partial topological dynamics}

In this final section we use the results of section 2 to characterize partial actions of a compact space by clopen sets whose associated partial skew group ring is simple. More specifically, we will prove the following theorem.

\begin{theorem}\label{t4} Let $\theta=(\{X_t\}_{t \in G}, \{h_t\}_{t \in G})$ be a partial action of a group $G$ on a compact space $X$ such that for each $t\in G$, $X_t$ is a clopen set. Then the partial skew group ring $\mathcal{C}(X)\rtimes_\alpha G $, where $\mathcal{C}(X)$ denotes the continuous complex-valued functions on $X$,  is simple if, and only if, $\theta$ is topologically free and minimal. 
\end{theorem}

\begin{obs} Partial actions on the Cantor set by clopen subsets are exactly the ones for which the enveloping space is Hausdorff (see \cite{EGG}). 
\end{obs}  

\begin{obs} %By \cite{Ex}, $\mathcal{C}(X)\rtimes_\alpha G $ is associative, and
Since the partial action acts on clopen sets, each $D_t$ has a unit. Hence, we can use Theorem \ref{t1} to prove the above theorem.
\end{obs}

Before we proceed, recall that there is a correspondence between partial actions on a locally compact Hausdorff space $X$ and partial actions on the C*-algebra of continuous complex-valued functions vanishing at infinity, $\mathcal{C}_0(X)$, (see \cite{ELQ, BG} for example). Namely, if $\theta=(\{X_t\}_{t \in G}, \{h_t\}_{t \in G})$ is a partial action on $X$, then $\alpha=(\{D_t\}_{t \in G}, \{\alpha_t\}_{t \in G})$, where $D_t=\mathcal{C}_0(X_t)$ and $\alpha_t(f):=f\circ h_{t^{-1}}$, is a partial action of $G$ on $\mathcal{C}_0(X)$. Simplicity of the associated C*-partial crossed product was studied in \cite{ELQ}, and a version of the above theorem for partial actions of abelian groups was given in \cite{DG}. Below we will recall the relevant definitions and make the proper adaptations of the ideas in \cite{DG} to the case at hand.

\begin{definition} A topological partial action $\theta=(\{X_t\}_{t \in G}, \{h_t\}_{t \in G})$ is topologically free if for all $t \neq 0$ the set $F_t=\{ x \in X_{t^{-1}}: h_t(x) = x \}$ has empty interior and is minimal if there is no proper, open invariant subset of $X$ ($U\subset X $ is invariant if $h_t(U\cap X_{t^{-1}} )\subseteq U$ for all $t\in G$).
\end{definition}

\begin{proposition} A partial action $\theta=(\{X_t\}_{t \in G}, \{h_t\}_{t \in G})$ on a compact space $X$ is minimal if, and only if, $\mathcal{C}(X)$ is $G-$simple.
\end{proposition}

\dem
The proof of this can be found in \cite{ELQ}.
\fim

\begin{proposition} Suppose that $\theta=(\{X_t\}_{t \in G}, \{h_t\}_{t \in G})$ is a topologically free partial action. Then $\mathcal{C}(X)\delta_0$ is maximal commutative in $\mathcal{C}(X)\rtimes_\alpha G $.
\end{proposition}

\dem
Suppose that $\mathcal{C}(X)\delta_0$ is not maximal commutative. Then there exists a non-zero function $f_t$ and $t\in G$, with $t\neq 0$, such that $f_t \delta_t\cdot f \delta_0 =  f \delta_0\cdot f_t \delta_t$ for all $f \in \mathcal{C}(X)$, which is equivalent to $\alpha_t(\alpha_{t^{-1}}(f_t)f) \delta_t = f f_t \delta_t$,  for all $f \in \mathcal{C}(X)$, which in turn is equivalent to 

\begin{equation}\label{eq41}
f_t(x) f(h_{t^{-1}}(x)) = f(x) f_t(x), 
\end{equation} for all $f \in \mathcal{C}(X)$ and $x \in X_t$.

Now, since $f_t$ is non-zero, there exists $x\in X_t$ such that $f_t(x) \neq 0$ and the continuity of $f_t$ implies that there exists an open set $U\subseteq X_t$ such that $f_t$ is non-zero in $U$. Since the partial action is topologically free there exists $y\in U$ such that $h_{t^{-1}}(y)\neq y$. Let $f\in \mathcal{C}(X)$ be such that $f(y)=1$ and $f(h_{t^{-1}}(y)) = 0$ (such a function exists by Urysohn's lemma). But then equation \eqref{eq41} above implies that $f_t(y)=0$, a contradiction.
\fim

\begin{proposition}
 If $\mathcal{C}(X) \rtimes_\alpha G$ is simple, then $\theta=(\{X_t\}_{t \in G}, \{h_t\}_{t \in G})$ is topologically free.
 \end{proposition}

\dem
The proof of this proposition is analogous to the proof of Proposition 4.7 in \cite{DG}.
\fim

\begin{obs} The three propositions above, combined with Theorem \ref{t1}, provide the proof of Theorem \ref{t4}.
\end{obs}

\end{document}